\documentclass{article}

\usepackage[utf8]{inputenc}
\usepackage[T1]{fontenc}
\usepackage{amsmath,amssymb}
\usepackage{indentfirst}
\usepackage{hyperref}

\newenvironment{theorem}{\begin{quote}\it}{\end{quote}}

\newcommand{\guil}[1]{\guillemotleft{}#1\guillemotright{}}

\begin{document}

\title{Graphs and matrices: A translation \\ of ``Graphok és matrixok'' by Dénes Kőnig (1931)\footnote{The talk ``Graphs and matrices'' was given at the seminar of the Lorand Eötvös Mathematics and Physics Society (Eötvös Loránd Matematikai és Fizikai Társulat) on March 26, 1931. The paper accompanying the talk originally appeared in the ``Matematikai és Fizikai Lapok'' (Mathematical and Physical Journal), volume 38, 1931. The original paper in Hungarian is available online at \url{http://real-j.mtak.hu/7307/}. This translation has been published with permission by the János Bolyai Mathematical Society (Bolyai János Matematikai Társulat). Thanks to Anna Gujgiczer and Naomi Arnold for providing feedback on the initial drafts of the translation.}}
\author{Gábor Szárnyas \\ \texttt{szarnyasg@gmail.com}}
\date{}

\maketitle

\begin{abstract}
    This paper, originally written in Hungarian by D\'{e}nes K\H{o}nig in 1931, proves that in a bipartite graph, the minimum vertex cover and the maximum matching have the same size. This statement is now known as K\H{o}nig's theorem. The paper also discusses the connection of graphs and matrices, then makes some observations about the combinatorial properties of the latter.
\end{abstract}

Let $G$ be a (finite) bipartite graph.
This means that all closed paths %\footnote{The original text used the term \emph{lines} (``vonalak'') to denote \emph{sequences of edges}.}
in $G$ are of even length, in other words, the vertices of $G$ can be partitioned into two sets $\Pi_1$ and $\Pi_2$ such that all edges in $G$ connect a vertex in $\Pi_1$ with a vertex in $\Pi_2$.
Let $M$ be the maximal number of the edges in $G$ which do not have a common vertex.%
\footnote{Such a set of edges is now called an \emph{independent edge set} or a \emph{matching}.}
If vertices $A_1, A_2, \ldots, A_v$ in $G$ are such that all edges in $G$ are incident to one of these vertices, we say that $A_1, A_2, \ldots, A_v$ \emph{cover} the edges of $G$.\footnote{Instead of \emph{cover} (``lefogják''/``lefedik''), the original text used the term \emph{exhaust} (``kimerítik'').}

We prove that \emph{the edges of $G$ can be covered with $M$ vertices}.

Let
$$ K = (P_1 Q_1, P_2 Q_2, \ldots, P_M Q_M) $$
be the set of $M$ edges such that, in accordance with the definition of $M$, the vertices $P_i, Q_i (M = 1, 2, \ldots, M)$ are distinct.
Let vertices $P_i$ belong to set $\Pi_1$ and vertices $Q_i$ to $\Pi_2$%(i.e.\ $P_i \in \Pi_1, Q_i \in \Pi_2$)
, and let
$$\Pi_1' = (P_1, P_2, \ldots, P_M), \qquad \Pi_2' = (Q_1, Q_2, \ldots, Q_M),$$
in a way that $\Pi_1'$ is a subset of $\Pi_1$ and $\Pi_2'$ is a subset of $\Pi_2$. We base our proof on the notion of a \guil{\emph{$K$-path}}.

A \emph{$K$-path} in $G$ is a path (which is open and does not have repeating vertices) $A_1 A_2 \ldots A_{2r}$ in which the second, fourth, \ldots, $2v$\textsuperscript{th}, \ldots, and penultimate edge, i.e.\ the edges $A_2 A_3, A_4 A_5, \ldots, A_{2v} A_{2v+1}, \ldots, A_{2r-2} A_{2r-1}$ all belong to $K$.
First, we prove the following lemma:

\begin{theorem}
    There is no \emph{$K$-path} in $G$ that connects a vertex in $\Pi_1 - \Pi_1'$ to another vertex in $\Pi_2 - \Pi_2'$.
\end{theorem}

Suppose $U$ would be such a path, then by removing the edges of $U \cap K$ from $K$ and adding the edges of $U \setminus K$ (where the size of the latter set is greater by 1), we would obtain $M + 1$ edges which do not share a vertex.
This contradicts the maximal nature of $M$.

Now we define a subset of $\Pi_1' + \Pi_2'$, $\Pi' = (R_1, R_2, \ldots, R_M)$:
let $\alpha$ be any of $1, 2, \ldots, M$ and let $R_\alpha = Q_\alpha$ if some $K$-path connects a vertex in $\Pi_1 - \Pi_1'$  with $Q_\alpha$; if there is no such $K$-path, let $R_\alpha = P_\alpha$.
This way, $\Pi'$ contains an endpoint of all the edges in $K$.
We prove that the set $\Pi'$ of $M$ vertices covers the edges of $G$, i.e.\ -- given that $PQ$ is an arbitrary edge of $G$ (where $P$ is in $\Pi_1$ and $Q$ is in $\Pi_2$) -- either $P$ or $Q$ is in $\Pi'$.
Our proof distinguishes between four cases:

\emph{Case 1.}
    Let $P$ belong to $\Pi_1 - \Pi_1'$ and $Q$ to $\Pi_2 - \Pi_2'$.
    By adding this $PQ$ edge to $K$, we would obtain $M + 1$ edges which do not share a vertex. This contradicts the maximal property of $M$, therefore, this case is not possible.

\emph{Case 2.}
    Let $P$ belong to $\Pi_1 - \Pi_1'$ and $Q$ to $\Pi_2'$.
    Then, $Q = Q_\alpha$, where $\alpha = 1, 2, \ldots, \text{ or } M$ and edge $PQ$ alone forms a $K$-path which connects vertex $P$ in $\Pi_1 - \Pi_1'$ to $Q = Q_\alpha$.
    Therefore, $Q = Q_\alpha$ belongs to $\Pi'$.

\emph{Case 3.}
    Let $P$ belong to $\Pi_1'$ and $Q$ to $\Pi_2 - \Pi_2'$.
    Then, $P = P_\alpha$, where $\alpha = 1, 2, \ldots, \text{ or } M$.
    If there were a $K$-path which connects some vertex $P_0$ of $\Pi_1 - \Pi_1'$ with $Q_\alpha$, then by adding edges $Q_\alpha P_\alpha$ and $P_\alpha Q$ we would derive a $K$-path, which connects $P_0$ with $Q$.
    However, this is impossible according to our lemma. 
    Therefore, there is no $K$-path which connects a vertex in $\Pi_1 - \Pi_1'$ to $Q_\alpha$.
    Therefore, $P = P_\alpha$ belongs to $\Pi'$.

\emph{Case 4.}
    Let $P$ belong to $\Pi_1'$ and $Q$ to $\Pi_2'$.
    Let e.g.\ $P = P_\alpha, Q = Q_\beta$. If $\alpha = \beta$, then trivially either $P = P_\alpha$ or $Q = Q_\alpha$ belong to $\Pi'$.
    Therefore, assume that $\alpha \neq \beta$.
    Either $P = P_\alpha$ belongs to $\Pi'$ or there exists a $K$-path, which connects some vertex $P_0$ of $\Pi_1 - \Pi_1'$ with $Q_\alpha$; in the latter case, by adding edges $Q_\alpha P_\alpha$ and $P_\alpha Q_\beta$ to this $K$-path, we derive a $K$-path which connects $P_0$ with $Q_\beta$ in way that $Q = Q_\beta$ belongs to $\Pi'$.

With this, we have indeed proved that if a bipartite graph maximally has $M$ edges which do not have a common vertex, then the edges of $G$ can be covered by $M$ vertices.
Therefore, if $m$ is the minimal number of vertices which cover the vertices of $G$, then $m \leqq M$.

It is obvious that the opposite is also true: $m \geqq M$.
If viz.\ $e_1, e_2, \ldots, e_M$ are edges which do not have a common vertex and vertices $A_1, A_2, \ldots, A_m$ cover the edges of the graph, then all of edges $e_1, e_2, \ldots, e_M$ end in one of the vertices $A_1, A_2, \ldots, A_m$; and as these do not have any common endpoints, indeed $m \geqq M$.

With this, we have proved that $m = M$.
To summarize, our main result can be stated as follows:

\begin{theorem}
    In a bipartite graph, the minimal number of vertices covering all edges is equal to the maximal number of edges which do not have a common vertex.
\end{theorem}

Turning our attention to the application of this theorem on matrices, let
$$\Vert a_{ik} \Vert \qquad (i = 1, 2, \ldots, p; k = 1, 2, \ldots, q)$$
be any matrix where regarding the value of a certain element, we only consider whether it ``vanishes'' or not.\footnote{The term \emph{vanishes} is a literal translation of the word ``eltűnik'' (meaning vanishes, disappears, or fades away). In recent works on sparse linear algebra, these elements are called \emph{zero elements}, while \emph{non-vanishing elements} are called \emph{non-zero elements}.}
This matrix corresponds to a bipartite graph as follows.
Each of the $p$ rows corresponds to one of the vertices $P_1, P_2, \ldots, P_p$, each of the $q$ columns corresponds to  one of vertices $Q_1, Q_2, \ldots, Q_q$; furthermore we create a $P_i Q_k$ edge iff the \emph{corresponding} $a_{ik}$ element does not vanish.
We do not create any other edges.
This way, we construct a bipartite graph $G$.

The notion that vertices cover the edges of $G$ clearly means that the set of vertices corresponding to these rows and columns (in general: \emph{lines}) contain all non-vanishing elements of the matrix.
Meanwhile, the notion that certain edges do not have a common vertex means that the elements corresponding to these edges do not lie on the same line.

Overall, our result for matrices can be summarized as follows:

\begin{theorem}
    For any matrix, the minimal number of lines which contain all non-vanishing elements is equal to the maximal number of non-vanishing elements which pairwise do not lie on the same line.
\end{theorem}

It is obvious that the term \guil{non-vanishing} can be substituted with any property of the elements, therefore this theorem expresses a purely combinatorial property of matrices (two-dimensional tables) where elements can be any objects (not just numbers).

Finally, we mention that our results are closely related to the research on determinants by \mbox{\textsc{Frobenius}} and on graphs by \mbox{\textsc{Menger}}. We will discuss these connections later.

\hfill \emph{Dénes Kőnig.}

\end{document}